\documentclass[12pt]{article}
\usepackage{amsmath,amsthm,amsfonts,amssymb}
\def\End{{\mathrm {End}}}

\def\id{{\mathrm {id}}}

\def\Vir{{\mathrm {Vir}}}

\def\epsilon{\varepsilon}

\def\la{\lambda}

\def\phi{\varphi}
\def\si{\sigma}

\newtheorem{theorem}{Theorem}[section]
\newtheorem{lemma}[theorem]{Lemma}

\newtheorem{proposition}[theorem]{Proposition}
\newtheorem{remark}[theorem]{Remark}

\def\Q{{\mathbb Q}}

\title{{\bf The classification of non-local chiral CFT with $c<1$}}

\author{
{\sc Yasuyuki Kawahigashi}\footnote{Supported in part by JSPS.}\\
Department of Mathematical Sciences\\
University of Tokyo, Komaba, Tokyo, 153-8914, Japan\\
e-mail: {\tt yasuyuki@ms.u-tokyo.ac.jp}\\
\vphantom{X}\\
{\sc Roberto Longo$^\dag$}\footnote{Supported in part by GNAMPA and
  MIUR.}\footnote{Supported in part by EU network ``Quantum Spaces -
  Noncommutative Geometry" HPRN-CT-2002-00280}\\
Dipartimento di Matematica\\
Universit\`a di Roma ``Tor Vergata''\\
Via della Ricerca Scientifica, 1, I-00133 Roma, Italy\\
e-mail: {\tt longo@mat.uniroma2.it}\\
\vphantom{X}\\
{\sc Ulrich Pennig$^\ddag$} and {\sc Karl-Henning Rehren$^\ddag$}\\
Institut f\"ur Theoretische Physik\\
Universit\"at G\"ottingen\\
37077 G\"ottingen, Germany\\
e-mail: {\tt pennig@theorie.physik.uni-goe.de}, \\
 {\tt rehren@theorie.physik.uni-goe.de}}

\begin{document}
\maketitle

\begin{center}\sl Dedicated to Hans-J\"urgen Borchers on the occasion of
his 80th birthday \end{center} 
\vskip10mm 

\begin{abstract}
All non-local but relatively local irreducible extensions of Virasoro
chiral CFTs with $c<1$ are classified. The classification, which is a
prerequisite for the classification of local $c<1$ boundary CFTs on a
two-dimensional half-space, turns out to be $1$ to $1$ with certain
pairs of $A$-$D$-$E$ graphs with distinguished vertices.
\end{abstract}

\vspace{10mm}

\section{Introduction}

Non-local chiral conformal quantum field theories have gained renewed
interest because they give rise to local CFT on the two-dimensional Minkowski
halfspace $x>0$ (boundary CFT, BCFT), and vice versa \cite{LR2}. 

More precisely, a BCFT contains {\em chiral} fields which generate
a net $A$ of local algebras on the circle, such that $A_+(O) =
A(I)\vee A(J)$ are the chiral BCFT observables localized in the double
cone $O=I\times J\equiv \{(t,x): t+x\in I, t-x\in J\}$ where $I > J$
are two open intervals of the real axis (= the pointed circle). The
two-dimensional {\em local} fields of the BCFT define a net of
inclusions $A_+(O)\subset B_+(O)$ subject to locality, conformal
covariance, and certain irreducibility requirements. 

If $A$ is assumed to be completely rational \cite{KLM}, then there is
a 1 to 1 correspondence \cite{LR2} between Haag dual BCFTs 
associated with a given chiral net $A$, and non-local chiral extensions 
$B$ of $A$ such that the net of inclusions $A(I)\subset B(I)$ is 
covariant, irreducible and relatively local, i.e., $A(I)$ commutes
with $B(J)$ if $I$ and $J$ are disjoint. The correspondence is given
by the simple relative commutant formula 
$$B_+(O) = B(K)'\cap B(L)$$
where $O=I\times J$ as before, $K$ is the open interval between $I$ and
$J$, and $L$ is the interval spanned by $I$ and $J$. Conversely,
$$B(L) = \bigvee_{I\subset L,\; J\subset L,\; I>J} B_+(O).$$

BCFTs which are not Haag dual are always intermediate between $A_+$
and a Haag dual BCFT.

The classification of local BCFTs on the two-dimensional halfspace is
thus reduced to the classification of non-local chiral extensions,
which in turn \cite{LR1} amounts to the classification of $Q$-systems
(Frobenius algebras) in the $C^*$ tensor category of the superselection
sectors \cite{DHR} of $A$. The chiral nets $A=\Vir_c$ defined by the
stress-energy tensor (Virasoro algebra) with $c<1$ are known to be
completely rational, so the classification program just outlined can
be performed. 

{\em Local} chiral extensions of $\Vir_c$ with $c<1$ have a
direct interpretation as local QFT models of their own. Their classification
has been achieved previously (\cite{KL1}, see Remark 2.3) by imposing
an additional condition \cite{LR1} on the $Q$-system involving the
braided structure (statistics \cite{DHR}) of the tensor category. Of
course, the present non-local classification contains the local one. 

As in \cite{KL1} we exploit the fact that the tensor subcategories
of the ``horizontal'' and of the ``vertical'' superselection sectors
of $\Vir_c$ with $c<1$ are isomorphic with the tensor categories of
the superselection sectors of $SU(2)$ current algebras. (The braiding
is different, however.) We therefore first classify the $Q$-systems in
the latter categories (Sect.\ 1), and then proceed from $Q$-systems in the
subcategories to $Q$-systems in the tensor categories of all sectors of 
$\Vir_c$ (Sect.\ 2). Thanks to a cohomological triviality result
\cite{KL2}, the classification problem simplifies considerably, and
essentially reduces to a combinatorial problem involving the Bratteli
diagrams associated with the local subfactors $A(I)\subset B(I)$,
combined with a ``numerological'' argument concerning Perron-Frobenius
eigenvalues.  

In the last section, we determine the vacuum Hilbert spaces of the
non-local extensions and of the associated BCFT's thus classified.

\section{Classification of irreducible non-local extensions of the
  $SU(2)_k$-nets} 
\label{SU2}

As an easy preliminary, we first classify all irreducible, possibly
non-local, extensions of the $SU(2)_k$-nets on the circle.  Consider
the representation category of the $SU(2)_k$-net and label the
irreducible DHR sectors as $\la_0,\la_1,\la_2,\dots,\la_k$ as usual.
(The vacuum sector is labeled as $\la_0$.)  Label this net as $A$
and an irreducible extension as $B$.  Since $A$ is completely rational
in the sense of \cite{KLM}, the index $[B:A]$ is automatically
finite by \cite[Proposition 2.3]{KL1}.  We need to classify
the irreducible $B$-$A$ sectors ${}_B \iota_A$, where $\iota$ is
the inclusion map.  Note that the $A$-$A$ sector
${}_A \bar\iota \iota_A$ gives
the dual canonical endomorphism of the inclusion and this
decomposes into a direct sum of $\la_j$ by \cite{LR1}.  Suppose we
have such a sector ${}_B \iota_A$, and consider the following
sequence of commuting squares.

$$\begin{array}{ccccccccc}
\End({}_A \id_A) & \subset & \End({}_A {\la_1}{}_A) & \subset &
\End({}_A {\la^2_1}{}_A) & \subset & \End({}_A {\la^3_1}{}_A) &
\subset & \cdots\\
\cap && \cap && \cap && \cap && \\
\End({}_B \iota_A) & \subset & \End({}_B \iota{\la_1}{}_A) & \subset &
\End({}_B \iota{\la^2_1}{}_A) & \subset & \End({}_B \iota{\la^3_1}{}_A) &
\subset & \cdots
\end{array}$$

The Bratteli diagram of the first row arises from reflections
of the Dynkin diagram $A_{k+1}$ as in usual subfactor theory, that is,
it looks like Fig. \ref{F1}, where each vertex is labeled with
irreducible sectors appearing in the irreducible decomposition of
${\la^k_1}$ for $k=0,1,2,\dots$.  (See \cite[Section 9.6]{EK}
for appearance of such graphs in subfactor theory.)

\unitlength 0.7mm
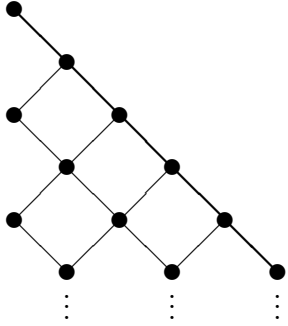
\begin{figure}[tb]
\begin{center}
\begin{picture}(90,70)
\thicklines
\put(10,60){\line(1,-1){50}}
\thinlines
\put(10,40){\line(1,-1){30}}
\put(10,20){\line(1,-1){10}}
\put(10,40){\line(1,1){10}}
\put(10,20){\line(1,1){20}}
\put(20,10){\line(1,1){20}}
\put(40,10){\line(1,1){10}}
\multiput(20,10)(20,0){3}{\circle*{3}}
\multiput(10,20)(20,0){3}{\circle*{3}}
\multiput(20,30)(20,0){2}{\circle*{3}}
\multiput(10,40)(20,0){2}{\circle*{3}}
\put(20,50){\circle*{3}}
\put(10,60){\circle*{3}}
\put(20,5){\makebox(0,0){$\vdots$}}
\put(40,5){\makebox(0,0){$\vdots$}}
\put(60,5){\makebox(0,0){$\vdots$}}
\end{picture}
\end{center}
\caption{The Bratteli diagram of the first row}
\label{F1}
\end{figure}

The Bratteli diagram of the second row also arises from reflections
of some Dynkin diagrams having the same Coxeter number as $A_{k+1}$
and starts with a single vertex because of the irreducibility
of ${}_B \iota_A$.  This gives a bipartite graph $G$, one of the
$A$-$D$-$E$ Dynkin diagram and its initial vertex $v$ as an invariant
of $\iota$, but 
the vertex $v$ is determined only up to the graph automorphism, 
so we denote the orbit of a vertex $v$ under such
automorphisms by $[v]$.  (Note that the Dynkin diagrams $A_n$,
$D_n$, and $E_6$ have non-trivial graph automorphisms of order 2.)
Also note that the graph $G$ is bipartite by definition.  A $B$-$A$
sector corresponding to an even vertex of $G$ might be equivalent
to another $B$-$A$ sector corresponding to an odd vertex of $G$.
Later it turns out that this case does not occur in the $SU(2)_k$-case,
but it does occur in the Virasoro case below.
Fig. \ref{F2} shows an example of $G$ and $v$ where $G$ is the
Dynkin diagram $E_6$.

\unitlength 0.7mm
\begin{figure}[tb]
\begin{center}
\begin{picture}(120,70)
\thinlines
\put(20,60){\line(1,-1){30}}
\put(10,50){\line(1,-1){40}}
\put(10,30){\line(1,-1){20}}
\put(70,40){\line(1,-1){10}}
\put(90,40){\line(1,-1){10}}
\put(10,50){\line(1,1){10}}
\put(10,30){\line(1,1){20}}
\put(10,10){\line(1,1){30}}
\put(30,10){\line(1,1){20}}
\put(80,30){\line(1,1){10}}
\put(100,30){\line(1,1){10}}
\put(30,10){\line(0,1){40}}
\put(90,30){\line(0,1){10}}
\multiput(10,10)(20,0){3}{\circle*{3}}
\multiput(20,20)(10,0){3}{\circle*{3}}
\multiput(10,30)(20,0){3}{\circle*{3}}
\multiput(20,40)(10,0){3}{\circle*{3}}
\multiput(70,40)(20,0){3}{\circle*{3}}
\multiput(80,30)(10,0){3}{\circle*{3}}
\put(10,50){\circle*{3}}
\put(30,50){\circle*{3}}
\put(20,60){\circle*{3}}
\put(10,5){\makebox(0,0){$\vdots$}}
\put(30,5){\makebox(0,0){$\vdots$}}
\put(50,5){\makebox(0,0){$\vdots$}}
\put(90,22){\makebox(0,0){$G$}}
\put(80,23){\makebox(0,0){$v$}}
\end{picture}
\end{center}
\caption{The Bratteli diagram of the second row}
\label{F2}
\end{figure}
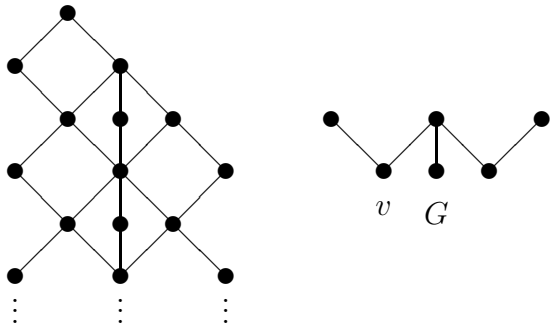

\begin{theorem}\label{SU2-class}
The pair $(G, [v])$ gives a complete invariant for irreducible
extensions of nets $SU(2)_k$, and an arbitrary pair $(G, [v])$ arises
as an invariant of some extension.
\end{theorem}

\begin{proof}
As in the proof of \cite[Proposition A.3]{BEK2}, we know that the
paragroup generated by $\iota$ is uniquely determined by $(G,[v])$ and
it is isomorphic to the paragroup of the Goodman-de la Harpe-Jones
subfactor given by $(G,[v])$, which was defined in 
\cite[Section 4.5]{GHJ}.  Then we obtain uniqueness of the $Q$-system
for the extension, up to unitary equivalence, as in
\cite[Theorem 5.3]{KL2}.  (We considered only a local extension
of $SU(2)_{28}$ corresponding to $E_8$ and its vertex having the
smallest Perron-Frobenius eigenvector entry there, but the same
method works for any $(G,[v])$.)

Any combination of $(G,[v])$ is possible as in the proof of
\cite[Lemma A.1]{BEK2}.  (We considered only the case of $E_7$ and 
its vertex having the smallest Perron-Frobenius eigenvector entry
there, but the same method works for any $(G,[v])$.
\end{proof}

\begin{remark}\label{Rem1}{\rm
Note that the pair $(G,[v])$ uniquely corresponds to an (isomorphism
class of) irreducible Goodman-de la Harpe-Jones subfactor.
So we may say that irreducible extensions of nets $SU(2)_k$ are
labeled with irreducible Goodman-de la Harpe-Jones subfactors
$N \subset M$ such that the inclusions
$A(I) \subset B(I)$ of the localized algebras are isomorphic
to $N \subset M$ tensored with a type III factor.
}\end{remark}

\section{Classification of non-local extensions of the Virasoro
nets $\Vir_c$ with $c<1$}

Let $A$ be the Virasoro net with central charge $c<1$.  As studied
in \cite[Section 3]{KL1}, it is a completely rational net.  We would
like to classify all, possibly non-local, irreducible extensions of
this net.  For $c=1-6/m(m+1)$, $m=3,4,5,\dots$, we label the
irreducible DHR sectors of the net $A$ as follows.  We have
$\si_{j,k}$, $j=0,1,\dots,m-2$,
$k=0,1,\dots,m-1$, with identification of
$\si_{j,k}=\si_{m-2-j,m-1-k}$.  (Our notation $\si_{j,k}$ corresponds
to $\la_{j+1,k+1}$ in \cite[Section 3]{KL1}.)  Our labeling gives that
the identity sector is $\si_{0,0}$ and the statistical dimensions
of $\si_{1,0}$ and $\si_{0,1}$ are $2\cos(\pi/m)$ and
$2\cos(\pi/(m+1))$, respectively.  We have $m(m-1)/2$ irreducible
DHR sectors.  We again need to classify the irreducible
$B$-$A$ sectors ${}_B \iota_A$, where $\iota$ is the inclusion map.  
Take such ${}_B \iota_A$ for a fixed $\Vir_c$ with $c=1-6/m(m+1)$
and we obtain an invariant as follows.

Consider the following sequence of commuting squares as in the
Section \ref{SU2}.
$$\begin{array}{ccccccccc}
\End({}_A \id_A) & \subset & \End({}_A {\si_{1,0}}{}_A) & \subset &
\End({}_A {\si_{1,0}^2}{}_A) & \subset & \End({}_A {\si_{1,0}^3}{}_A) &
\subset & \cdots\\
\cap && \cap && \cap && \cap && \\
\End({}_B \iota_A) & \subset & \End({}_B \iota{\si_{1,0}}{}_A) & \subset &
\End({}_B \iota{\si_{1,0}^2}{}_A) & \subset
& \End({}_B \iota{\si_{1,0}^3}{}_A) &
\subset & \cdots
\end{array}$$

 From the Bratteli diagram of the second row, we obtain a graph $G_1$
and its vertex $v_1$ as in Section \ref{SU2}.  The graph $G_1$
is one of the $A$-$D$-$E$ Dynkin diagrams and has the Coxeter
number $m$.  We also use $\si_{0,1}$ instead of $\si_{1,0}$
in this procedure and obtain a graph $G_2$ and its vertex $v_2$.
The graph $G_2$ is one of the $A$-$D$-$E$ Dynkin diagrams and has
the Coxeter number $m+1$. The quadruple $(G_1, [v_1], G_2, [v_2])$
is an invariant for $\iota$.  (The notation $[\cdot]$ means the
orbit under the graph automorphisms as in Section \ref{SU2}.)
Note that one of the graphs $G_1, G_2$ must be of type $A$ 
because the $D$ and $E$ diagrams have even
Coxeter numbers. We then prove the following classification theorem.

\begin{theorem}\label{Vir-class}
The quadruple $(G_1, [v_1], G_2, [v_2])$ gives a complete invariant
for irreducible extensions of nets $\Vir_c$, and an arbitrary quadruple,
subject to the conditions on the Coxeter numbers as above, arises
as an invariant of some extension.
\end{theorem}

We will distinguish certain sectors by their dimensions.
For this purpose, we need the following technical lemma on
the values of dimensions, which we prove before
the proof of the above theorem.

\begin{lemma}\label{trig}
Let $m$ be a positive odd integer and $G$ one of the $A$-$D$-$E$
Dynkin diagrams having a Coxeter number $n$ with $|n-m|=1$.
Take a Perron-Frobenius eigenvector $(\mu_a)_a$ for the graph $G$,
where $a$ denotes a vertex of $G$.  
Set $d_j=\sin(j \pi/m)/\sin(\pi/m)$ for $j=1, 2, \dots, m-1$.
Then the sets
$$\{\mu_a/\mu_b\mid a, b {\rm\ are\ vertices\ of\ }G\}$$
and
$$\{d_2,d_3,\dots,d_{m-2}\}$$
are disjoint.
\end{lemma}

\begin{proof}
If $m=1,3$, then the latter set is empty, so we may assume $m\ge5$.
Note that the value $1$ is not in the latter set.

Suppose a number $\omega$ is in the intersection and
we will derive a contradiction.  Then $\omega$
is in the intersection of the cyclotomic fields
$\Q(\exp(\pi i/m))$ and $\Q(\exp(\pi i/n))$, which is $\Q$ since
$(2m,2n)=2$ and $\Q(\exp(2\pi i/2))=\Q$.
Suppose $d_j$ is equal to this $\omega$.  We may and do assume
$2\le j \le (m-1)/2$.  We have
$$\omega=\frac{\zeta^j-\zeta^{-j}}{\zeta-\zeta^{-1}}
=\zeta^{j-1}+\zeta^{j-3}+\zeta^{j-5}+\cdots+\zeta^{-j+1},$$
where $\zeta=\exp(2\pi i/(2m))$. 

First assume that $j$ is even.  We note $(m-2, 2m)=1$ since $m$
is odd.  Then the map $\sigma:\zeta^k \mapsto \zeta^{k(m-2)}$ 
for $k=0,1,\dots, 2m-1$ gives an element of the Galois group
for the cyclotomic extension $\Q\subset\Q(\zeta)$.  We have
$$\sigma(\omega)=
\zeta^{(j-1)(m-2)}+\zeta^{(j-3)(m-2)}+\zeta^{(j-5)(m-2)}+
\cdots+\zeta^{(-j+1)(m-2)}.$$
Here the set
$$\{\zeta^{(j-1)(m-2)}, \zeta^{(j-3)(m-2)}, \zeta^{(j-5)(m-2)},
\dots,\zeta^{(-j+1)(m-2)}\}$$
has $j$ distinct roots of unity containing $\zeta^{m-2}$
and it is a subset of
$$Z=\{\zeta^k\mid k=1,3,5,\dots, 2m-1\}.$$
The set
$$\{\zeta^{j-1},\zeta^{j-3},\zeta^{j-5},\dots,\zeta^{-j+1}\}$$
is the unique subset having $j$ distinct elements of $Z$
that attains the maximum of
${\rm Re}\; \sum_{k=1}^j \alpha_k$
among all subsets
$\{\alpha_1,\alpha_2,\dots,\alpha_j\}$ having
$j$ distinct elements of $Z$.
However, we have $m>3$, which implies $j\le (m-1)/2< m-2$, thus
the complex number $\zeta^{m-2}$ is not in the above unique set, and
thus the sum
$$\zeta^{(j-1)(m-2)}+\zeta^{(j-3)(m-2)}+\zeta^{(j-5)(m-2)}+
\cdots+\zeta^{(-j+1)(m-2)}$$ cannot be equal to
$$\zeta^{j-1}+\zeta^{j-3}+\zeta^{j-5}+\cdots+\zeta^{-j+1},$$
which shows that $\omega$ is not fixed by $\sigma$, so $\omega$
is not an element of $\Q$, which is a contradiction.

Next we assume that $j$ is odd.  We now have that
$$\zeta^{2(j-1)/2}+\zeta^{2(j-3)/2}+\cdots+\zeta^{2(1-j)/2}\in\Q.$$
Since $(m,(m-1)/2)=1$, the map $\sigma:\zeta^{2k} \mapsto \zeta^{k(m-1)}$ 
for $k=0,1,\dots, m-1$ gives an element of the Galois group
for the cyclotomic extension $\Q\subset\Q(\zeta^2)$. Since
$m-1>j-1$, $\sigma(\omega)$ contains a term $\zeta^{m-1}$ which does not
appear in $\omega$.  Then by an argument similar to the above case of
even $j$, we obtain a contradiction.
\end{proof}

We now start the proof of Theorem \ref{Vir-class}.

\begin{proof}
Without loss of generality, we may assume that $m$ is odd and hence
that the graph $G_1$ is $A_{m-1}$.
(Otherwise, the graph $G_2$ is $A_m$, and we can switch
the symmetric roles of $G_1$ and $G_2$.)  Note that a sector
corresponding to an even vertex of $A_{m-1}$ can be equivalent
to another sector corresponding to an odd vertex of $A_{m-1}$.

The tensor category
having the irreducible objects $\{\si_{0,0},\si_{1,0},\dots,
\si_{m-2,0}\}$ is isomorphic to the representation category of
$SU(2)_{m-2}$, thus each of the irreducible objects is labeled with
a vertex of the Dynkin diagram $A_{m-1}$.  Let $\si_{j,0}$ be one
of the two sectors corresponding to $[v_1]$.  We choose $j$ to be
even, and then $j$ is uniquely determined.  Let $\Delta$ be the
set of the irreducible $B$-$A$ sectors arising from the 
decomposition of ${}_B \iota\si_{1,0}^{2k}{}_A$ for all $k$.
Note that $\Delta$ is a subset of the vertices of $G_1$.
Let $\tilde\iota$ be one of the $B$-$A$ sectors in $\Delta$
having the smallest dimension.
By the Perron-Frobenius theory and the definition of the graph
$G_1$, which is now $A_{m-1}$, we know that such $\tilde\iota$
is uniquely determined and that the set
$$\{d(\lambda)/d(\tilde\iota)\mid \la\in\Delta\}$$
is equal to 
$$\{\sin(k\pi/m)/\sin(\pi/m)\mid k=1,2,\dots,(m-1)/2\}.$$
The situation is illustrated in Fig. \ref{F3} where we also have
the graph $G$ which will be defined below.  The vertices corresponding
to the elements in $\Delta$ are represented as larger circles.

\unitlength 0.7mm
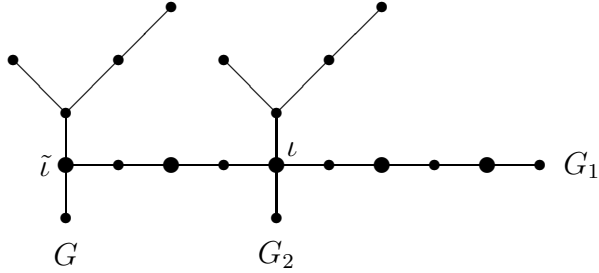
\begin{figure}[tb]
\begin{center}
\begin{picture}(120,60)
\thinlines
\put(10,40){\line(1,-1){10}}
\put(50,40){\line(1,-1){10}}
\put(20,30){\line(1,1){20}}
\put(60,30){\line(1,1){20}}
\put(20,10){\line(0,1){20}}
\put(60,10){\line(0,1){20}}
\put(20,20){\line(1,0){90}}
\multiput(20,20)(20,0){5}{\circle*{3}}
\multiput(30,20)(20,0){5}{\circle*{2}}
\multiput(20,10)(0,20){2}{\circle*{2}}
\multiput(60,10)(0,20){2}{\circle*{2}}
\multiput(10,40)(20,0){4}{\circle*{2}}
\multiput(40,50)(40,0){2}{\circle*{2}}
\put(20,3){\makebox(0,0){$G$}}
\put(60,3){\makebox(0,0){$G_2$}}
\put(118,20){\makebox(0,0){$G_1$}}
\put(63,23){\makebox(0,0){$\iota$}}
\put(16,20){\makebox(0,0){$\tilde\iota$}}
\end{picture}
\end{center}
\caption{The graphs $G_1, G_2, G$}
\label{F3}
\end{figure}

Now we consider the Bratteli diagram for
$\End(\tilde\iota\si_{0,1}^n)$ for $n=0,1,2,\dots$, and obtain a graph
$G$ and an orbit $[v]$ whose reflection gives this Bratteli diagram. 
The graph $G$ is one of the $A$-$D$-$E$ Dynkin diagrams and its
Coxeter number differs from $m$ by $1$.

We want to show that the pairs $(G_2,[v_2])$ and $(G,[v])$ are equal
as follows.
We first claim that the irreducible decomposition of
$\overline{\tilde\iota}\tilde\iota$ contains only sectors
among $\si_{0,k}$, $k=0,1,\dots,m-1$. Suppose that
$\overline{\tilde\iota}\tilde\iota$ contains $\si_{2l,k}$
with $l>0$ on the contrary.   By the Frobenius reciprocity, 
we have 
$$0< \langle \overline{\tilde\iota}\tilde\iota, \si_{2l,k}\rangle
=\langle \tilde\iota \si_{0,k}, \tilde\iota \si_{2l,0}\rangle.$$
By the above description of the graph $G_1$, we know that
$\tilde\iota \si_{2l,0}$ is irreducible and distinct from
$\tilde\iota$.   The assumption that $\tilde\iota\si_{2l,0}$ appears
in the irreducible decomposition of $\tilde\iota\si_{0,k}$ means that
the graphs $G$ and $G_1$ have a common vertex other than $\tilde\iota$
and this is impossible by Lemma \ref{trig}.  We have thus proved that 
the irreducible decomposition of
$\overline{\tilde\iota}\tilde\iota$ contains only sectors
among $\si_{0,k}$, $k=0,1,\dots,m-1$.

We know $\tilde\iota \si_{j,0}$ and $\iota$ are equivalent sectors since
they both are irreducible. To show $(G_2,[v_2])=(G,[v])$, we therefore
need to compare the irreducible decompositions of $\tilde\iota
\si_{j,0}\si_{0,1}^k$ and $\tilde\iota \si_{0,1}^k$ for
$k=0,1,2,\dots$.  Suppose that $\la$ is an irreducible sector
appearing in the decomposition of $\tilde\iota \si_{0,1}^k$ for some $k$.  
We have 
$$\langle \la \si_{j,0}, \la\si_{j,0}\rangle=
\langle \bar\la \la, \si_{j,0}^2\rangle.$$
Now the decomposition of $\bar\la\la$ 
contains only sectors among $\si_{0,l}$, $l=0,1,\dots,m-1$
as above.  Thus the only irreducible sector appearing in
decompositions of both $\bar\la\la$ and $\si_{j,0}^2$ is
the identity sector, which appears exactly once in the both.
We conclude that $\la \si_{j,0}$ is also irreducible.  Thus the
irreducible decompositions of $\tilde\iota \si_{j,0}\si_{0,1}^k$
and $\tilde\iota \si_{0,1}^k$ for $k=0,1,2,\dots$ are described
by the same Bratteli diagram and 
$(G,[v])$ and $(G_2,[v_2])$ are equal.

Then ${}_A \overline{\tilde\iota}\tilde\iota{}_A$ decomposes into
some irreducible sectors among $\si_{0,0},\si_{0,1},\dots,
\si_{0,m-1}$.  Since the tensor category having 
the irreducible objects
$$\{\si_{0,0}, \si_{0,1},\dots, \si_{0,m-1}\}$$
is isomorphic to the representation category of
$SU(2)_{m-1}$, we obtain uniqueness of $\tilde\iota$ for a given
$(G,[v])$ hence $(G_2,[v_2])$.  Then $\iota=\tilde\iota\si_{j,0}$ 
determines a $Q$-system uniquely, up to unitary equivalence.

We next prove a realization of a given $(G_1, [v_1], G_2, [v_2])$.
We continue to assume that $G_1$ is $A_{m-1}$.  Let $\si_{j,0}$ be one
of the two sectors corresponding to $[v_1]$ as above.
Using the tensor category
having the irreducible objects $\{\si_{0,0},\si_{0,1},\dots,
\si_{0,m-1}\}$, we have $\tilde\iota$ corresponding to $(G_2, [v_2])$
as in the proof of Theorem \ref{SU2-class}.  Set
$\iota=\tilde\iota\si_{j,0}$.  Then one can verify that this $\iota$
produces the quadruple $(G_1, [v_1], G_2, [v_2])$ by the
same argument as in the above one showing $G=G_2$.
\end{proof}

\begin{remark}{\rm 
By \cite[Theorem 4.1]{KL1}, we already know that the local extensions
among the above classification are labeled with
$(A_{n-1},A_n)$, $(A_{4n},D_{2n+2})$, $(D_{2n+2},A_{4n+2})$,
$(A_{10},E_6)$, $(E_6, A_{12})$,
$(A_{28},E_8)$ and $(E_8, A_{30})$ for $(G_1, G_2)$ and
the vertices $v_1, v_2$ are those having the smallest
Perron-Frobenius eigenvector entries.
}\end{remark}

\begin{remark}{\rm 
As in Remark \ref{Rem1}, we may say that
irreducible extensions of the Virasoro nets with $c<1$ are
labeled with pairs of irreducible Goodman-de la Harpe-Jones subfactors
having the Coxeter numbers differing by $1$.
}\end{remark}

\begin{remark}{\rm 
The graphs $G_1$ and $G_2$ are by definition bipartite, thus excluding
the tadpole graphs which also have Frobenius norm $<2$. Tadpole diagrams
arise by pairwise identification of the vertices of $A_m$ diagrams when
$m=2n$ is even. Indeed, when $\iota:A\hookrightarrow B$ equals
$\sigma_{n-1,j}:A\to A$, the even vertices of $G_2=A_{2n}$ pairwise
coincide as $B$-$A$ sectors with the odd vertices, so that the fusion
graph for multiplication by $\si_{0,1}$ is $T_n$. The invariant $G_2$
in these cases is $A_m$, nevertheless.  

As an example, consider the case $m=4$, that is,
$c=7/10$. In this case, we have six irreducible DHR sectors for the
net $\Vir_{7/10}$. The graphs $G_1$ and $G_2$ are automatically $A_3$
and $A_4$, respectively, so we have four possibilities for the
invariant $(G_1, [v_1], G_2, [v_2])$.  If $(G_1, [v_1], G_2, [v_2])$ is
as in case (1) of Fig. \ref{F4}, then the sector $\iota$ is given
by $\si_{1,1}$, thus the four vertices of the graph $A_4$ give only
two mutually inequivalent $B$-$A$ sectors. That is, the fusion graph
of the $B$-$A$ sectors for multiplication by $\si_{0,1}$ is the
tadpole graph as in Fig. \ref{F5}. 

In case (2) of Fig. \ref{F4} we have four mutually inequivalent
$B$-$A$ sectors for the graph $G_2=A_4$ for the sector $\iota$ given
by $\si_{0,1}$, and the fusion graph is also $A_4$.
}\end{remark}

\unitlength 0.9mm
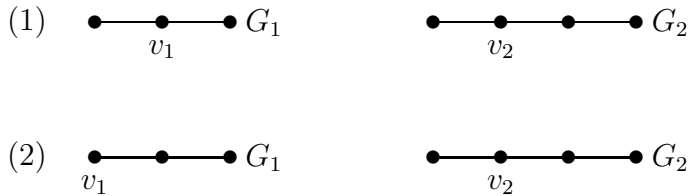
\begin{figure}[tb]
\begin{center}
\begin{picture}(100,40)
\thinlines
\put(10,30){\line(1,0){20}}
\put(10,10){\line(1,0){20}}
\put(60,30){\line(1,0){30}}
\put(60,10){\line(1,0){30}}
\multiput(10,10)(10,0){3}{\circle*{2}}
\multiput(10,30)(10,0){3}{\circle*{2}}
\multiput(60,10)(10,0){4}{\circle*{2}}
\multiput(60,30)(10,0){4}{\circle*{2}}
\put(35,10){\makebox(0,0){$G_1$}}
\put(35,30){\makebox(0,0){$G_1$}}
\put(95,10){\makebox(0,0){$G_2$}}
\put(95,30){\makebox(0,0){$G_2$}}
\put(0,10){\makebox(0,0){$(2)$}}
\put(0,30){\makebox(0,0){$(1)$}}
\put(20,26){\makebox(0,0){$v_1$}}
\put(70,26){\makebox(0,0){$v_2$}}
\put(10,6){\makebox(0,0){$v_1$}}
\put(70,6){\makebox(0,0){$v_2$}}
\end{picture}
\end{center}
\caption{Two cases of $(G_1, [v_1], G_2, [v_2])$ for $m=4$, $c=7/10$}
\label{F4}
\end{figure}

\unitlength 0.9mm
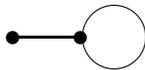
\begin{figure}[tb]
\begin{center}
\begin{picture}(40,30)
\thinlines
\put(10,10){\line(1,0){10}}
\put(25,10){\circle{10}}
\multiput(10,10)(10,0){2}{\circle*{2}}
\end{picture}
\end{center}
\caption{The tadpole graph $T_2$}
\label{F5}
\end{figure}

\medskip

\section{The canonical endomorphism}

We want to determine, viewed as a representation of the subtheory 
$A={\rm Vir}_c$, the vacuum Hilbert space of the local boundary 
conformal QFT associated with each of the non-local extensions $B$, 
classified in the previous section. This representation is given by a 
DHR endomorphism $\theta$ of $A$ whose restriction to a local algebra 
$A(I)$ (where $\theta$ is localized in the interval $I$) coincides 
with the canonical endomorphism $\bar\iota\iota$ of the subfactors 
$A(I)\subset B(I)$ classified above. We are therefore interested in the
computation of $\bar\iota\iota$. 

By the equality of local and global intertwiners, and by reciprocity,
the multiplicity of each irreducible DHR sector $\sigma$ within
$\theta$ equals the multiplicity of $\iota$ within $\iota\sigma$. We
therefore need to control the decomposition of $\iota\sigma$ into
irreducibles (``fusion'') for all DHR sectors. Because every
irreducible sector is a product $\sigma_{j,k} = \sigma_{j,0}\sigma_{0,k}$, 
and $\sigma_{j,0}$ and $\sigma_{0,k}$ are obtained from the
generators $\sigma_{1,0}$ and $\sigma_{0,1}$ by the recursion
$\sigma_{0,k+1} = \sigma_{0,k}\sigma_{0,1} \ominus \sigma_{0,k-1}$ and
likewise for $\sigma_{j+1,0}$, it suffices to control the fusion with
the generators.   

We know from the preceding section that the fusion of $\iota$ with the
generators $\sigma_{1,0}$ and $\sigma_{0,1}$ separately can be
described in terms of the two bi-partite graphs $G_1$ and $G_2$ such
that the vertices of the graphs represent irreducible $B$-$A$-sectors
and two vertices are linked if the corresponding sectors are connected
by the generator. $\iota$ corresponds to a distinguished vertex in
both graphs. Moreover, we have seen that the fusion of $\iota$ with
both generators can be described by the ``product graph'' $G=G_1\times
G_2$ with vertices $\lambda=(v_1\in G_1,v_2\in G_2)$ and
``horizontal'' edges linking $(v_1,v_2)$ with $(v_1',v_2)$ if $v_1$
and $v_1'$ are linked in $G_1$, and likewise for ``vertical'' edges
according to the graph $G_2$. Again, the vertices $\lambda$ represent
irreducible $B$-$A$-sectors, and $\iota$ is a distinguished vertex
of the product graph.  See Fig.\ \ref{F6} for an example.

\unitlength 1.0mm
\begin{figure}[tb]
\begin{center}
\begin{picture}(120,60)
\thinlines
\put(14,39){\line(1,0){50}}
\multiput(27,36)(10,0){5}{\line(1,0){6}}
\multiput(20,10)(0,10){3}{\line(1,0){50}}

\multiput(20,10)(10,0){6}{\line(0,1){20}}
\multiput(20,30)(10,0){6}{\line(-2,3){6}}
\multiput(20,30)(10,0){6}{\line(2,3){4}}

\multiput(20,10)(10,0){6}{\circle*{2}}
\multiput(20,20)(10,0){6}{\circle*{2}}
\multiput(20,30)(10,0){6}{\circle*{2}}
\multiput(24,36)(10,0){6}{\circle*{2}}
\multiput(14,39)(10,0){6}{\circle*{2}}

\put(108,20){\makebox(0,0){$G=G_1\times G_2$}}
\put(63,23){\makebox(0,0){$\iota$}}
\end{picture}
\end{center}
\caption{The graph $G_1\times G_2$ with $G_1=A_6$ and
  $G_2=D_5$. Different vertices may represent the same
  $B$-$A$-sector. $\iota$ may be any vertex of $G$.}
\label{F6}
\end{figure}
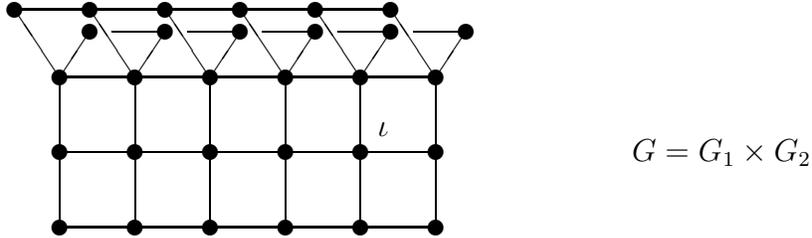

From the product graph $G$, the fusion of each of its vertices with any DHR 
sector can be computed in terms of vertices of $G$, i.e.,
$\lambda\sigma$ can be decomposed into irreducibles represented by the
vertices of $G$. But different vertices may represent identical 
$B$-$A$-sectors; we only know that within each horizontal or vertical 
subgraph, the even vertices represent pairwise inequivalent sectors,
and so do the odd vertices, cf.\ Remark 3.5. In order 
to compute the canonical endomorphism, we have to determine all 
identifications between vertices of $G$ as $B$-$A$-sectors.

We continue to assume that $m$ is odd and hence $G_1 = A_{m-1}$. 
The Coxeter number $m'=2n$ of $G_2$ is either $m+1$ or $m-1$. We 
exploit the fact that $\sigma_{m-2,0}$ and $\sigma_{0,m'-2}$ represent the 
same DHR sector $\tau$, and that $\tau$ is simple (it has dimension 1). 
Hence fusion with $\tau$, as a horizontal sector $\sigma_{m-2,0}$,
yields an automorphism $\alpha_1$ of the graph $G_1$ such that
$\alpha_1(v_1)= v_1\sigma_{m-2,0}$, and, as a vertical sector,
similarly yields an automorphism $\alpha_2$ of $G_2$. It follows that
the vertices $\lambda=(v_1,v_2)$ and $\alpha(\lambda) =
(\alpha_1(v_1),\alpha_2(v_2))$ of the product graph represent the same
sectors.  

Because $\tau$ connects even vertices of $G_1$ with odd ones,
$\alpha_1$ must be the unique non-trivial automorphism of $G_1$. 
To determine $\alpha_2$, one may use the above-mentioned recursion to 
compute the fusion of the vertices of $G_2$ with $\tau=\sigma_{0,m'-2}$.
We find that $\alpha_2$ is the unique non-trivial automorphism if $G_2$ is
either an $A$ graph or $E_6$ or $D_{2n+1}$, and it is trivial if $G_2$
is $E_7$, $E_8$, or $D_{2n}$.

We now claim that the identifications due to $\alpha$ give {\em all}\/
pairs of vertices of $G$ which represent the same $B$-$A$-sector.

\begin{proposition}
The graph
$$(G_1\times G_2)/(\alpha_1\times\alpha_2)$$
is the fusion graph of $\iota$ with respect to $\sigma_{1,0}$ and
$\sigma_{0,1}$, i.e., its vertices represent inequivalent irreducible
$B$-$A$-sectors, and its horizontal and vertical edges correspond to
fusion with the two generators.    
\end{proposition}

\begin{proof} 
By the Perron-Frobenius theory, the dimensions of the
$B$-$A$-sectors represented by the vertices $\lambda=(v_1,v_2)$ of $G$ are
common multiples of $\nu(v_1)\mu(v_2)$ where $\nu(v_1)$ and $\mu(v_2)$
are the components of the Perron-Frobenius eigenvectors $\nu$ of $G_1$ and
$\mu$ of $G_2$. Let now $\lambda=(v_1,v_2)$ and
$\lambda'=(v_1',v_2')$ be two vertices of $G$ which represent the same
$B$-$A$-sector. Then clearly 
$$\nu(v_1)\mu(v_2)=\nu(v_1')\mu(v_2').$$
If $v_1$ it at distance $j$ from an extremal vertex $\tilde v_1$ of
$G_1$, then $\tilde \lambda=(\tilde v_1,v_2)$ is a subsector of
$\lambda\sigma_{j,0}$ and consequently of $\lambda'\sigma_{j,0}$. Hence
$\tilde \lambda$ is equivalent to some subsector $(\tilde v_1',v_2')$
of $\lambda'\sigma_{j,0}$, implying $\mu(v_2)/\mu(v_2') = \nu(\tilde
v_1')/\nu(\tilde v_1) = d_k$ for some $k$.   
Lemma 3.2 tells us that this is only possible if $d_k=1$. It follows
that $\mu(v_2')=\mu(v_2)$ and $\nu(v_1')=\nu(v_1)$. 

This means in particular that $v_1'=v_1$ or $v_1'=\alpha_1(v_1)$, and
that $v_2$ and $v_2'$ and $\alpha_2(v_2)$ are all even or all odd. If 
$v_1'=v_1$, then $\lambda$ and $\lambda'$ are two even or two odd vertices within
the same vertical subgraph representing the 
same sector. This is only possible if $\lambda=\lambda'$. If on the other hand 
$v_1'=\alpha_1(v_1)$, then the same argument applies to $\alpha(\lambda)$ 
and $\lambda'$, giving $\lambda'=\alpha(\lambda)$.
\end{proof}

Having determined the fusion graph, it is now straightforward to
compute (as described above) the canonical endomorphism for every
possible position of $\iota$ as a distinguished vertex of the fusion
graph, and hence to determine the vacuum Hilbert space for each local
boundary conformal QFT with $c<1$. 

We display below the canonical endomorphism $\theta_{\tilde v_1,\tilde v_2}$
whenever $\tilde v_1$ and $\tilde v_2$ are extremal vertices of $G_1$ and
$G_2$. All other cases are then easily obtained by the following argument: 
If $v_1$ is at distance $j$ from an extremal vertex $\tilde v_1$ of $G_1$,
then $v_1 = \tilde v_1\sigma_{j,0}$. If $v_2$ is at distance $k$
from the extremal vertex $\tilde v_2$ on the same ``leg'' of $G_2$, then
$v_2 = \tilde v_2\sigma_{0,k}$. (If $v_2$ is the trivalent vertex
of the $D$ or $E$ graphs, then this is true for each of the three
legs.) It then follows that $(v_1,v_2) = (\tilde v_1,\tilde v_2)
\sigma_{j,k}$, and hence 
$$\theta_{v_1,v_2} = \theta_{\tilde v_1,\tilde v_2}\sigma_{j,k}^2.$$ 

The canonical endomorphisms $\theta_{\tilde v_1,\tilde v_2}$ for all pairs
of extremal vertices of $G_1$ and $G_2$ are listed in the following
table. 
\bigskip
\bigskip

\begin{tabular}{l|c|c||l}
$G_2$ & $m'$ & {\rm dist.} & $\theta_{\tilde v_1,\tilde v_2}$ \cr
\hline
$A_n$ & $n+1$ & $-$ & $\sigma_{0,0}$ \cr
$D_n$ & $2n-2$ & 1 & $\sigma_{0,0} \oplus \sigma_{0,4} \oplus
\sigma_{0,8} \oplus \ldots \oplus \sigma_{0,4[n/2]-4}$ \cr
$D_n$ & $2n-2$ & $n-3$ & $\sigma_{0,0} \oplus \sigma_{0,2n-4}$ \cr
$E_6$ & $12$ & $1$ & $\sigma_{0,0} \oplus \sigma_{0,4} \oplus \sigma_{0,6} 
\oplus \sigma_{0,10}$ \cr
$E_6$ & $12$ & $2$ & $\sigma_{0,0} \oplus \sigma_{0,6}$\cr
$E_7$ & $18$ & $1$ & $\sigma_{0,0}\oplus\sigma_{0,4}\oplus
\sigma_{0,6}\oplus\sigma_{0,8}\oplus\sigma_{0,10}
\oplus\sigma_{0,12}\oplus\sigma_{0,16}$ \cr
$E_7$ & $18$ & $2$ & $\sigma_{0,0}\oplus\sigma_{0,6} \oplus\sigma_{0,10}
\oplus\sigma_{0,16}$ \cr
$E_7$ & $18$ & $3$ & $\sigma_{0,0}\oplus\sigma_{0,8} \oplus \sigma_{0,16}$ \cr
$E_8$ & $30$ & $1$ & \hskip-1.75mm
$\begin{array}{r}\sigma_{0,0}\oplus\sigma_{0,4}\oplus
\sigma_{0,6}\oplus\sigma_{0,8}\oplus 2\sigma_{0,10}\oplus
\sigma_{0,12}\oplus 2\sigma_{0,14} \oplus \\ \oplus
\sigma_{0,16}\oplus 2\sigma_{0,18}\oplus\sigma_{0,20}\oplus
\sigma_{0,22}\oplus\sigma_{0,24}\oplus
\sigma_{0,28}\end{array}$ \cr
$E_8$ & $30$ & $2$ &
$\sigma_{0,0}\oplus\sigma_{0,6}\oplus\sigma_{0,10}\oplus\sigma_{0,12}\oplus
\sigma_{0,16}\oplus\sigma_{0,18}\oplus\sigma_{0,22} \oplus \sigma_{0,28}$ \cr
$E_8$ & $30$ & $4$ & $\sigma_{0,0} \oplus \sigma_{0,10} \oplus \sigma_{0,18}
\oplus \sigma_{0,28}$ \cr
\end{tabular}
\bigskip

{\bf Table 4.1.} The canonical endomorphisms $\theta_{\tilde v_1,\tilde v_2}$ 
for all pairs of extremal vertices of $G_1$ and $G_2$. The entry in
the third column indicates the distance of $\tilde v_2$ from the
trivalent vertex, i.e., the length of the ``leg'' of $G_2$ on which
$\tilde v_2$ is the extremal vertex.
\bigskip
\medskip

The local chiral extensions classified earlier \cite{KL1} are
precisely those cases where $G_2$ is $A$, $D_{2n}$, $E_6$, or $E_8$,
and both $v_1$ and $v_2$ are extremal vertices (on the respective
longest leg in the $D$ and $E$ cases). 

In the non-local cases, the local algebras of the associated BCFT on
the half-space are the relative commutants as described in the
introduction. Note that, in order to determine the resulting
factorizing chiral charge structure \cite{LR2} of the local fields,
more detailed information about the DHR category and the $Q$-system is
needed, than the simple combinatorial data exploited in this work.

\bigskip
\bigskip

\noindent{\bf Acknowledgments.}
A part of this work was done during visits of the first-named
author to Universit\`a di Roma ``Tor Vergata'' and Universit\"at
G\"ottingen, and he thanks for their hospitality.  The authors
gratefully acknowledge the financial support of GNAMPA-INDAM and MIUR
(Italy), EU network ``Quantum Spaces - Noncommutative Geometry'', and
Grants-in-Aid for Scientific Research, JSPS (Japan). 

\end{document}